\documentclass{article}
\usepackage{graphicx,upquote}
\begin{document}

\markboth{J. L. Aurentz and L. N. Trefethen}{Chopping a Chebyshev series}

\title{Chopping a Chebyshev series}
\author{JARED L. AURENTZ and
LLOYD N. TREFETHEN}

\maketitle

\begin{abstract}
Chebfun and related software projects
for numerical computing with functions
are based on the idea that
at each step of a computation, a function $f(x)$ defined
on an interval $[a,b]$ is ``rounded'' to a prescribed
precision by constructing a Chebyshev series and chopping it
at an appropriate point.  Designing a chopping algorithm with
the right properties proves to be a
surprisingly complex and interesting problem.  We describe
the chopping algorithm introduced in Chebfun Version 5.3 in 2015
after many years of discussion and the considerations that
led to this design.
\end{abstract}

\section{Introduction: the construction problem}
Floating point arithmetic is based on the idea that at each
step of a computation, a number is rounded to a prescribed precision,
which in standard IEEE arithmetic is 53 bits
or about 16 digits~\cite{overton}.  By rounding at every step in this fashion,
one eliminates the combinatorial explosion of the lengths
of numerators and denominators that would 
occur in exact rational arithmetic: all numbers are represented
to the same relative accuracy and require the same storage.
It is ultimately for this reason that
virtually all computational science is carried out
in floating point arithmetic.

Chebfun and other related software systems that have
arisen in the past fifteen years are based on implementing
an analogous principle for functions as opposed to
numbers~\cite{battlestref,chebbook,cacm}.
If $f$ is a Lipschitz continuous function on $[-1,1]$, it has
a uniformly and absolutely convergent Chebyshev series
\begin{equation}
f(x) = \sum_{k=0}^\infty a_k T_k(x), \quad
a_k = {2\over \pi} \int_{-1}^1 {f(x) T_k(x)\over \sqrt{1-x^2}}\kern 1pt dx,
\label{series}
\end{equation}
where $T_k$ denotes the degree $k$ Chebyshev polynomial,
with $2/\pi$ changed to $1/\pi$ in the case of the coefficient
$a_0$~\cite{atap}.  If $f$ has several derivatives of smoothness or
better, the series typically converges to 16-digit accuracy after
tens or hundreds of terms, making computation with functions
in this form eminently practical when carried out with
stable and efficient algorithms, including solution of problems
such as rootfinding, integration, differentiation, 
minimization, and maximization. Piecewise-smooth functions can
be represented by concatenating such representations, with
intervals and subintervals $[a,b]$ other than $[-1,1]$ handled by the obvious
linear change of variables applied to (\ref{series}).

\begin{figure}
\begin{center}
\includegraphics[width=3.43in]{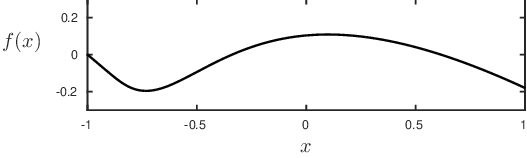}~\kern 13pt\\[-5pt]
\kern 13pt~\includegraphics[width=3.15in]{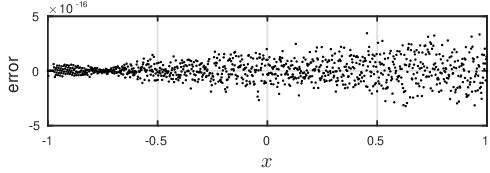}
\end{center}
\vskip -1.5em
\caption{\small\label{functionf}Above, the function $f$ of (\ref{examplefun})
represented in Chebfun by a polynomial of degree $165$.
Below, the error at\/ $1000$ equally spaced points in $[-1,1]$.}
\end{figure}

For example, consider the function
\begin{equation}
f = 3e^{-1/(x+1)} - (x+1), \quad x\in [-1,1],
\label{examplefun}
\end{equation}
which is $C^\infty$ but not analytic.
In Chebfun, in
$0.05$ secs.\ on a typical 2015 desktop computer, the following
commands plot $f$ and determine its zeros and maximum value in
$[-1,1]$.  The plot appears in Figure~\ref{functionf}.
\begin{verbatim}
>> f = chebfun(@(x) 3*exp(-1./(x+1))-(x+1));
>> plot(f)
>> roots(f)
   ans =
     -1.000000000000000
     -0.338683188672833
      0.615348950784159
>> max(f)
   ans = 0.108671573231256
\end{verbatim}

\begin{figure}
\begin{center}
\includegraphics[width=3.25in]{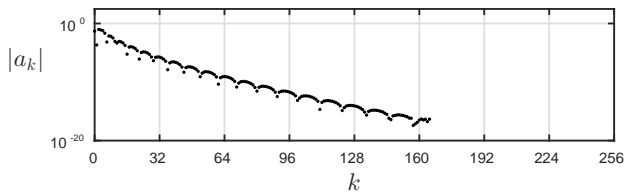}
\end{center}
\vskip -1.5em
\caption{\small\label{functionfcoeffs}Chebyshev
coefficients $|a_k|$ for the Chebfun representation of $f$.}
\end{figure}

\noindent Figure~\ref{functionf} also shows the error in
the Chebfun approximation to this function at $1000$ equally
spaced points in $[-1,1]$.  Evidently the approximation has about
15 digits of accuracy relative to the global scale of $f$.
Chebfun represents $f$ by a polynomial of degree 165, which
corresponds to a polynomial interpolant on a 166-point grid:
\begin{verbatim}
>> length(f)
   ans = 166
\end{verbatim}
This function representation is called a {\em chebfun} (with a
lower-case {\em c}).
Figure~\ref{functionfcoeffs} shows the absolute values of
the Chebyshev coefficients $a_k$ on a log scale.
This paper is about how representations like this are determined, not
only when a chebfun like $f$ is constructed from scratch, but as
further operations on chebfuns are carried out such as 
$f+g$, $fg$, $\exp(f)$, $1/f$, and so on.  (For the last of these,
poles are introduced at the zeros of $f$ and the
new function is represented by a chebfun with three smooth pieces, each
of which is called a {\em fun}.)
The same algorithm is applied for solution of ordinary differential
equation boundary-value problems (ODE BVPs) via Chebyshev
spectral collocation (Section~8).

\begin{figure}[p]
\begin{center}
\includegraphics[width=3in]{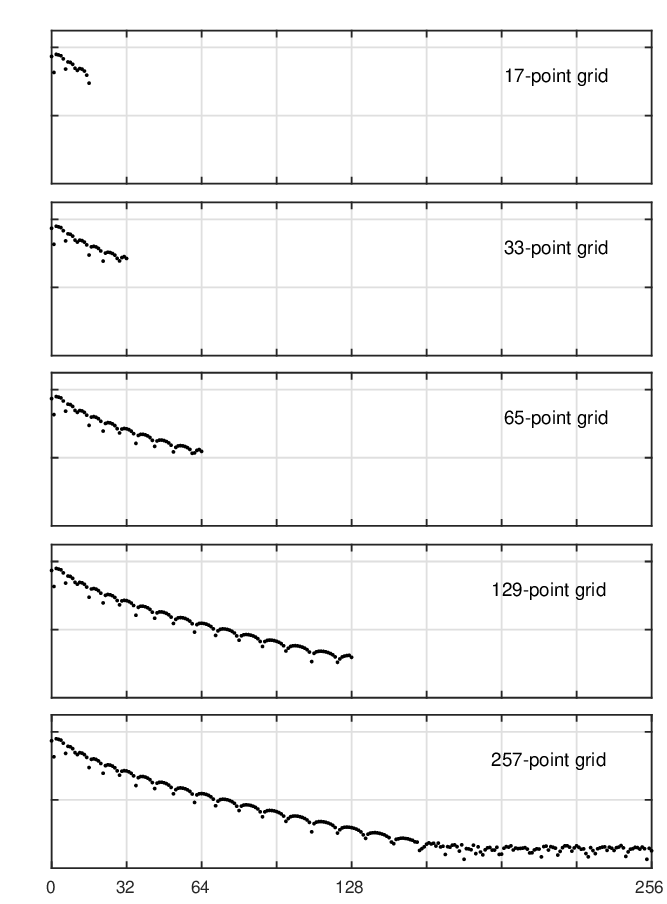}
\end{center}
\vskip -1em
\caption{\small\label{seriesfig}Chebfun construction process for the function $f$
of\/ $(\ref{examplefun})$.  The grid is refined repeatedly
until a plateau of
rounding errors near the level of machine
precision is detected, and then the series is chopped at
the beginning of the plateau to give the result
shown in Figure~$\ref{functionfcoeffs}$.}
\end{figure}

The core of the Chebfun construction process is
sampling a function on successively finer grids. 
In its standard mode of operation,
the constructor samples a function
on Chebyshev grids with $17, 33, 65,\dots$ points, computing
the Chebyshev coefficients of the polynomials of degrees
$16, 32, 64, \dots$ in sequence that interpolate the data on
each grid.  (The $(n+1)$-point
Chebyshev grid consists of the points
$\cos(\kern .7pt j\pi/n),$ $0\le j \le n$.)
When coefficients fall to the level of machine precision,
the refinement stops and the series is chopped.
Figure~\ref{seriesfig} shows how this process
plays out for this function $f$.
Grids of 17, 33, 65, and 129 points are tried first, and then
on the grid of 257 points, the Chebyshev coefficients reach
a plateau at a level
close to machine precision.
This corresponds to a Chebyshev series of
degreee 256, which the constructor then chops to degree 165.

\begin{figure}[p]
\begin{center}
\includegraphics[scale=.62]{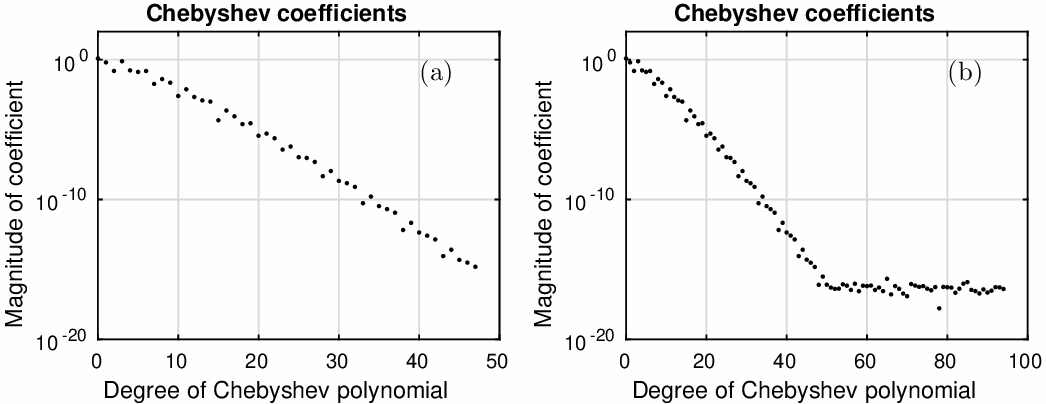}
\end{center}
\vskip -1em
\caption{\label{doublelength}\small The {\tt'doublelength'} flag constructs a chebfun
of twice the adaptively determined degree.  Combining this with
{\tt plotcoeffs} provides a simple way for users to
see how Chebfun has chosen to represent a function.
(a): {\tt plotcoeffs(chebfun(ff))}.
(b): {\tt plotcoeffs(chebfun(ff,'doublelength'))}.}
\end{figure}

A simple way for a user to see where Chebfun has chopped a series
is to call {\tt plotcoeffs} after constructing a function both
in the usual way and with the \verb|'doublelength'| flag, which
constructs it with twice the adaptively selected degree.  For example,
Figure~\ref{doublelength} plots such results for the function
$f(x) = \exp(\sin(\pi x))$ based on the commands
\begin{verbatim}
>> ff = @(x) exp(sin(pi*x));
>> f = chebfun(ff); plotcoeffs(f,'.')
>> f2 = chebfun(ff,'doublelength'); plotcoeffs(f2,'.')
\end{verbatim}

In brief, then, this paper is about the design of an automated procedure
for getting from Figure~\ref{seriesfig} to Figure~\ref{functionfcoeffs}.  It is based on more
than a decade of experience during which millions of chebfuns
have been constructed by ourselves and other
Chebfun users around the world.  This
experience has shown that there are surprisingly many challenges
along the way to making this straightforward-seeming construction
process fully reliable.

\section{Simplification: coefficients only}
To begin, we mention a simplifying assumption that applies throughout
our work.  This is that
chopping decisions are made entirely on the basis of Chebyshev
coefficients, not function values.  Loosely speaking, we have
chosen to regard the task
of getting from Figure~\ref{seriesfig} to Figure~\ref{functionfcoeffs}
as the ``whole problem'' of chebfun construction.

This construction in coefficient space is not
the only way in which the problem might have been approached.  
In principle, one might decide how fine a grid is needed 
by constructing an approximation on one grid and then evaluating it at
off-grid points.
Mathematically, this could lead to interesting questions associated
with pointwise accuracy of approximations and expansions.
On the other hand, many off-grid points would have to be
tested for reliability.  We would not claim that constructing
functions via coefficients is the only reasonable approach,
but for better or worse, this is the route Chebfun has
followed from the beginning.
One of the benefits of
this approach is that the human eye, looking at plots like
Figures~\ref{functionfcoeffs} and~\ref{seriesfig}, can help
judge how well an algorithm is working---and
Section~5 of this paper is built around such plots.

There is one place where the principle of
coefficients-only construction is broken.
After a candidate chebfun is
constructed, to confirm that it is valid,
Chebfun executes a code called {\tt sampleTest}
that samples the chebfun at two arbitrary points in the
interval of definition to make sure it matches
the function being constructed.  The idea here is to catch
the case in which the function being sampled contains an isolated high-order
Chebyshev polynomial component that aliases to a low-order one
on Chebyshev grids.  For example, the Chebyshev
polynomial $T_{128}(x)$ takes the value $1$ at every point of
the Chebyshev grids with 17, 33, and 65 points, so a pure
coefficient-based construction process would confuse this
function with the constant $1$.  For this function,
Chebfun gets the right answer since
{\tt sampleTest} duly fails for the candidate chebfuns
constructed on each of the grids just mentioned, until
finally the correct chebfun is identified on
the grid of 129 points.
\begin{verbatim}
>> length(chebfun(@(x) cos(128*acos(x))))
   ans = 129
\end{verbatim}
Of course, in principle a function might arise 
that outwits {\tt sampleTest}, but we have never seen this happen.
The possibility of being fooled by special data has been familiar
to numerical analysts since the 1960s, when de Boor,
Lyness and others
pointed out that adaptive quadrature routines can always
be fooled by an integrand that happens to take, say, the value~0
at every sample point~\cite{deboor}.\footnote{Any deterministic
adversary could be outwitted, with probability ${\approx}\kern .7pt 1$, by 
the use of randomized tests.  However, Chebfun
avoids randomization, because it is too disturbing to
users if results are occasionally unrepeatable, even if the
differences involve unimportant
choices at the level of machine precision.}

\section{Relative scales and convergence tolerance}
Following the principle
of floating-point arithmetic, all algorithmic decisions in
Chebfun are based on relative quantities.\footnote{The
solution of ODE IVPs is an exception, for our
algorithms are based on {\tt ode113} and other Matlab codes that
are not scale-invariant, as described
in Section~8.  We are considering ways to make them
scale-invariant nonetheless.}
Thus, for example,
$2^{500} f$ and $2^{-500} f$ are represented by polynomials
of exactly the same degree 165 as $f$ itself:
\begin{verbatim}
>> length(chebfun('2.^(500)*(3*exp(-1./(x+1))-(x+1))'))
   ans = 166
>> length(chebfun('2.^(-500)*(3*exp(-1./(x+1))-(x+1))'))
   ans = 166
\end{verbatim}
A further check of coefficients confirms that they
are exactly the same in all three cases apart from
factors of $2^{500}$.
To achieve scale-invariance, Chebfun makes decisions about
chopping a sequence of coefficients
relative to the largest coefficient in absolute value.
Scale-related irregularities will not
arise unless a function nears the overflow or underflow limits 
around $2^{\pm 1024}$.

A different question is what relative accuracy to seek in constructing
functions numerically.  Chebfun's design makes use of
a user-adjustable {\tt eps} parameter that determines
convergence, which is set by default to machine 
epsilon, $2^{-52}$.  For most Chebfun calculations, we have long
recommended that users leave {\tt eps} at this value rather than increasing
it in the hope of a speedup.\footnote{{\tt eps} can be
increased by the command
{\tt chebfunpref.setDefaults('eps',1e-8)} and returned to
its default value with
{\tt chebfunpref.setDefaults('factory')}.}
The reason is that 
Chebfun usually manages to exploit the piecewise smoothness of
functions, so that halving the accuracy requirement from 16 to 8 digits,
say, would typically at best halve the total cost.
We are motivated by the analogy with floating-point arithmetic,
where the consensus of decades of scientific computing is
that it is usually best to work in a fixed relatively fine precision
rather than try to make savings by fine-tuning the
precision.\footnote{For enthusiastic advocacy
of tuned precision, see~\cite{palmer}.}

There are two contexts in which it may most often
be desirable to loosen
{\tt eps} in Chebfun calculations.  One is in solving
differential equations, where the work may scale worse
than linearly as the grid is refined in search of greater accuracy---and,
moreover,
matrices may arise whose ill-conditioning forces the loss of
some digits in any case.  As we shall describe in
Section~8, Chebfun by default loosens {\tt eps} to about
$10^{-12}$ for ODE BVPs and IVPs.
The other is in working with noisy data.
If one is constructing functions from data only accurate to six digits,
say, then it will probably be important to set {\tt eps}
to a number on the order of $10^{-6}$.

\begin{figure}
\begin{center}
\includegraphics[width=3.4in]{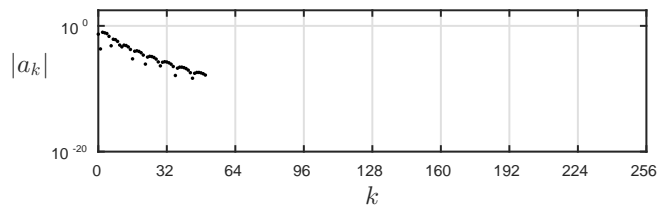}
\end{center}
\vskip -1.5em
\caption{\small\label{fig1e8}Chebyshev coefficients $|a_k|$ for the Chebfun representation
of $f$ constructed with ${\tt eps} = 10^{-8}$.  Compare
Figure~$\ref{functionfcoeffs}$, based on the default precision.}
\end{figure}

For example, here is a construction of a chebfun for the function
$f$ of (\ref{examplefun}) with precision $10^{-6}$ instead of
the default value.  The degree of the chebfun reduces
from $165$ to $50$, and $\max(f)$ is now accurate to 7 digits.
Figure~\ref{fig1e8} shows the Chebyshev coefficients.
A plot of the function itself
looks indistinguishable from Figure~\ref{functionf}a.
\begin{verbatim}
>> f = chebfun(@(x) 3*exp(-1./(x+1))-(x+1), 'eps', 1e-6);
>> length(f)
   ans = 51
>> max(f)
   ans = 0.108671567726459
\end{verbatim}
Further examples of loosening of {\tt eps} are presented
in Section~5 (Figure~\ref{fig9}).

\section{standardChop algorithm}
We now present the algorithm implemented in the Chebfun code
{\tt standardChop}, introduced in Version 5.3
in November 2015.  Some of aims in the
design of this algorithm have been these:
\begin{enumerate}
\item
Be simple enough to explain to people
(and ourselves!).\footnote{Folkmar Bornemann evoked
an embarrased laugh at the
Chebfun and Beyond conference in 2012 when he projected on
the screen a comment he had found at the heart of the
Version~4 Chebfun constructor Matlab code: {\tt \% Why do we do this?} .}
\item 
Be simple enough to reason about mathematically.
\item
In particular, have as few input and output parameters as possible.
\item
Unify diverse algorithms previously applied in various corners of
Chebfun (ODE BVPs, ODE IVPs, ``trigfuns'' for periodic
functions, 2D functions$,\dots$).
\item
Maintain accuracy as close as possible to 16 digits, or 10 digits for
ODEs (by default).
\item
Nevertheless, detect plateaus of rounding errors at a slightly higher
level if necessary and chop accordingly.
\item 
At the same time, get an extra digit or two when this is cheap.
\item
Adjust in a systematic way to user-specified looser tolerances,
e.g.\ for computing with functions contaminated by noise.
\item
Don't examine too many more coefficients than are eventually kept.
\item 
Get the right answer when the function really is a low-degree polynomial.
\item 
Avoid anomalies where Chebfun makes an ``obviously wrong'' chop.
\end{enumerate}
For full details about {\tt standardChop}, see the 
code listing (which includes careful comments)
in the appendix.  We do not claim that
this algorithm is optimal, merely that it is a good candidate in
a complex situation.  Indeed, the larger purpose of this article
is not so much to advocate a particular algorithm as to delineate some
of the surprisingly many considerations that arise in trying
to generalize the notion of floating point arithmetic from
numbers to functions.  In the remaining sections, the
motivation for the
design of {\tt standardChop} will be further explained via
many examples.

{\tt standardChop} takes two inputs, a real or
complex number sequence {\tt coeffs} of
length $n\ge 17$, 
and a relative tolerance ${\tt tol} \in (0,1)$, which is
normally set equal to Chebfun's general tolerance {\tt eps}.
It produces one output: an integer ${\tt cutoff}\in [1,n]$.\footnote{In
a project like Chebfun
with 100,000 lines of code and two dozen contributors, temptations
to make things complicated are ever-present.  Holding the numbers of
inputs and outputs to two and one, respectively, was an act of will.}
If ${\tt cutoff} \le n-1$, the constructor is ``happy''
and {\tt coeffs} should be chopped to length {\tt cutoff}.\ \ If
${\tt cutoff} = n$, the constructor is ``unhappy'' and
a longer coefficient sequence is needed.

Note that Matlab's convention of beginning indexing at 1 rather
than~0 is a potential source of confusion.
The sequence {\tt coeffs} is
indexed from 1 to $n$, so the description below is framed
in those terms.
In the standard application this will correspond to
a Chebyshev series with coefficients from degree $0$ to $n-1$,
and the final index {\tt cutoff} to be retained will
correspond to degree ${\tt cutoff}-1$.

{\tt standardChop} proceeds in three steps.

\medskip

{\em Step\/ $1$.  Compute the upper envelope of\/ {\tt coeffs} and normalize.}
The input sequence {\tt coeffs}
is replaced by the nonnegative, monotonically nonincreasing sequence 
${\tt envelope}_j = \max_{j\le k \le n} |{\tt coeffs}_k|$. (Note that
the use of the absolute value makes the chopping algorithm
applicable to complex functions as well as real ones.)  If
${\tt envelope_1}\ne 0$, this sequence is then normalized
by division by ${\tt envelope}_1$ to give a 
nonnegative, monotonically nonincreasing sequence whose first element
is $1$.  The output of {\tt standardChop} will depend
only on {\tt envelope}, so this first step entails a substantive algorithmic
decision: to assess a sequence only on its rate of
decay, ignoring any oscillations along the way.

\medskip

{\em Step\/ $2$.  Search for a plateau.}
The algorithm now searches for a sufficiently long, sufficiently
flat plateau of sufficiently small coefficients.  If
no such plateau is found, the construction process is
unhappy: {\tt cutoff} is set to $n$ and the algorithm
terminates.  A plateau can be as high as
${\tt tol}^{2/3}$ if it is perfectly flat but need not be
flat at all if it is as low as ${\tt tol}$.
Precisely, a plateau is defined as a stretch
of coefficients ${\tt envelope}_j,\dots, {\tt envelope}_{j_2}$
with $j\ge 2$ and $j_2 = \hbox{round}(1.25 j + 5)\le n$
with the property
$$
{{\tt envelope}_{j_2}\over {\tt envelope}_j} \ge r  =
3 \left(1 - {\log({\tt envelope}_j)\over \log({\tt tol})}\right).
$$
The integer {\tt plateauPoint} is set to $j-1$, where $j$ is the
first point that satisfies these conditions.

\medskip

{\em Step\/ $3$.  Chop the sequence near the beginning of the plateau.}
Having identified an index {\tt plateauPoint} that is followed by
a plateau, one might think that the code would simply set
${\tt cutoff} = {\tt plateauPoint}$ and terminate.  Such a procedure
works well most of the time.  However, exploring its application
to hundreds of functions reveals some examples where 
{\tt plateauPoint} doesn't catch the true ``elbow'' of
the {\tt envelope} curve.  Sometimes,
it is clear to the eye that a little more accuracy could be achieved
at low cost
by extending the sequence a little further (algorithmic aim~7 in the
list of the last section).
Other times, if a plateau is detected just below the
highest allowed level ${\tt tol}^{2/3}$, it is clear to the eye
that the plateau actually begins at an earlier point slightly
higher than ${\tt tol}^{2/3}$.  To adjust for these cases, Step~3
sets {\tt cutoff} not simply to
{\tt plateauPoint}, but to the index just before the lowest point of the
{\tt envelope} curve as measured against a straight line
on a log scale titled downward with a slope
corresponding to a decrease over the range
by the factor ${\tt tol}^{1/3}$.
Figure~\ref{stepsfig} gives the idea,
and for precise details, see the code listing in the appendix.
\begin{figure}[t!]
\begin{center}
\includegraphics[scale=.85]{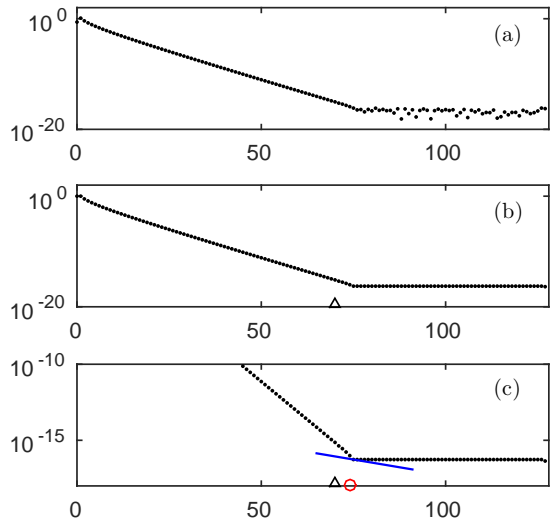}
\end{center}
\vskip -1.5em
\caption{\small\label{stepsfig}Sketch of the Chebfun construction
process for $f(x) = \log(1.1-x)$.
(a) After rejecting the $17$-, $33$-, and
$65$-point Chebyshev grids,
Chebfun computes coefficients on the $129$-point grid.
(b) In Step 1 of {\tt standardChop}, the
monotonically nonincreasing
normalized envelope of the coefficients is constructed, and
the plateau is found to be long, low, and level enough for chopping.
Step 2 picks ${\tt plateauPoint}=71$ as
the last point before the plateau, marked by a triangle at
position 70 on the axis (since the corresponding degree
is~70).
(c)
Step 3 finds the lowest coefficient relative to a line tilted
slightly downward, giving ${\tt cutoff} = 75$, marked by
a circle at position 74.\ \ For
this function the net effect of extending the series through
$|a_{74}|$ rather than $|a_{70}|$
is an improvement in accuracy by about one bit.} 
\end{figure}

\section{Examples and discussion}
We now present four figures each containing four
plots, for a total of sixteen examples selected to
illustrate various issues of chebfun construction.
We present these in the context of chebfuns constructed from
scratch, as in the example \verb|f = chebfun(@(x) exp(x)./(1+x.^2))|, but 
approximately the same results would arise in
computation with functions as in
\verb|x = chebfun(@(x) x)|, \verb|f = exp(x)./(1+x.^2)|, as discussed
in Section~8.
In each case we use the default interval $[-1,1]$.

\begin{figure}[t]
\begin{center}
\includegraphics[width=4.5in]{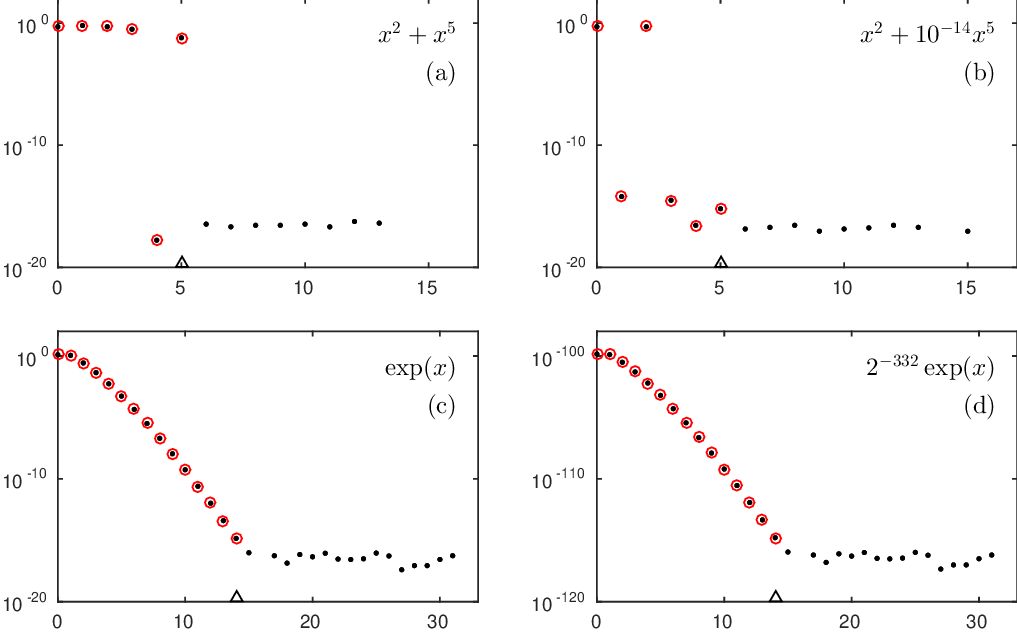}
\caption{\label{fig6}\small(a) Treatment of a polynomial.  (b) Capturing a small-magnitude
component.  (c) A non-polynomial function.  (d) Scale-invariance.}
\end{center}
\vskip .2in
\end{figure}

Each example is presented in the same format, as two sequences of
Chebyshev coefficients $|a_k|$ superimposed
on a log scale.  The open circles
show the coefficients of the chebfun as finally constructed.
The solid dots show the coefficients from the grid
of $2^K+1$ points on which the constructor has decided it is ``happy''.
Thus the sequence of solid dots is chopped to
obtain the final sequence of open circles.\footnote{Until Version~5.3,
the circles would not have matched the dots exactly, because a
sequence of coefficients was not simply chopped but aliased back
to the lower-order coefficients.  Aliasing improves accuracy at grid
points while worsening it between grid points, and in 
V5.3 we decided to drop it.}

The {\tt standardChop} indices {\tt cutoff} and {\tt plateauPoint} 
can be identified in these figures.  The final circle marks
${\tt cutoff}-1$, the highest degree of the chebfun finally
constructed.  The triangle plotted on the bottom axis marks
${\tt plateauPoint}-1$, the point just before the beginning
of the plateau that was detected.  Often these two points are
the same, indicating that Step 3 of {\tt standardChop} has had
no effect.  In several cases ${\tt cutoff}$ is slightly greater
than ${\tt plateauPoint}$, indicating that Step 3 has detected
that a little more accuracy can be attained at low cost, as
in Figure~\ref{stepsfig} (algorithmic
aim~7 in the list of Section~4).  In two cases, ${\tt cutoff}$
is much smaller than ${\tt plateauPoint}$, indicating that,
given that the series is to be chopped, it is more cost-effective
to chop it earlier (algorithmic aim~11).

\begin{figure}[t]
\begin{center}
\includegraphics[width=4.5in]{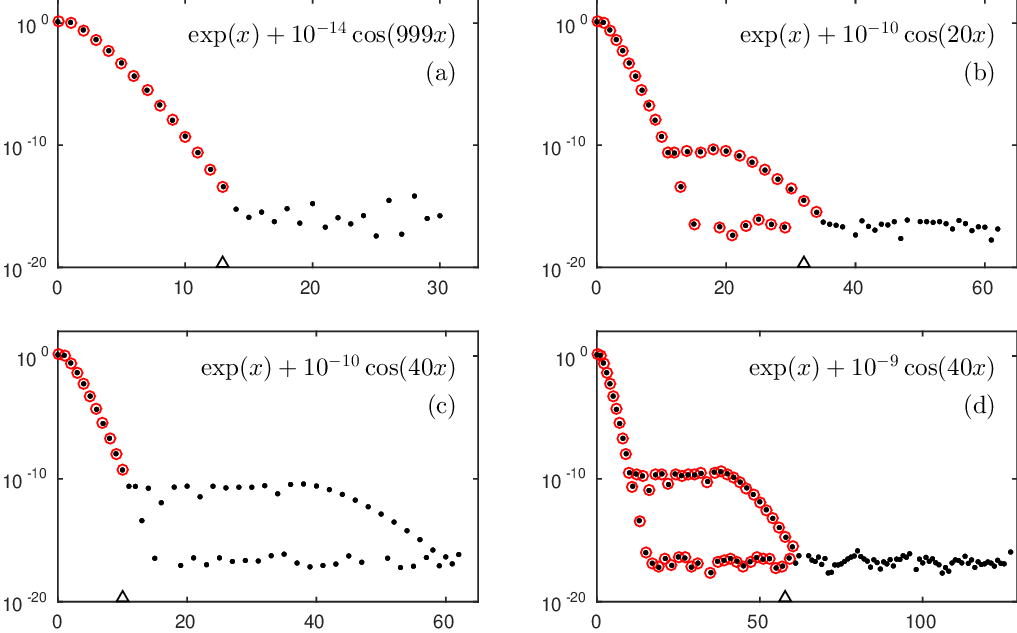}
\caption{\label{fig7}\small(a) A little noise leads to earlier chopping.
(b) This short stretch of small coefficients
is not treated as a plateau.
(c) This longer stretch is treated as a plateau.
(d) With higher amplitude, again it is not a plateau.}
\end{center}
\end{figure}

The mention of the term ``cost-effective'' highlights the fact that
throughout this design process,
choices of strategies and parameters have been made
that are somewhat arbitrary.  
How many coefficients that ``look
like a noise plateau'' should be sampled before the system decides
to look no further?  If coefficients have fallen close
to the ideal tolerance but it appears far more would be needed to
reach it, where should the series be chopped?
Can the eye be trusted if it judges that a chopping point
is ``obviously wrong''?
{\tt standardChop} embodies answers to such questions, and
the examples of this section
have been chosen to illustrate how these choices
play out in practice.  We hope we have struck a good balance
of accuracy, efficiency, and avoidance of disturbing anomalies.

We begin with examples \ref{fig6}a and \ref{fig6}b, selected to
remind us that Chebfun must deal effectively
with functions that are actually low-degree polynomials (algorithmic aim~10).
Usually there will be a high cliff between the coefficients
to be retained and those to be discarded, as with $x^2 + x^5$.
Sometimes, however, the cliff will be not so high, as
with $x^2 + 10^{-14} x^5$, and here the constructor must make a decision.
If $10^{-14}$ is changed to $10^{-16}$ in Figure \ref{fig6}b, it makes 
the other decision and chops the series after degree $2$.
(With the in-between value $10^{-15}$, not shown, it keeps the
$T_3$ coefficient but discards the slightly smaller $T_5$ one.)

Examples \ref{fig6}c and \ref{fig6}d begin the exploration of the
more usual case, functions that are not polynomials but will be
approximated by polynomials.  Example \ref{fig6}c is the most basic function
of all, $e^x$, and Example \ref{fig6}d divides this
by $2^{332}\approx 10^{100}$ to illustrate exact binary
scale-invariance.

\begin{figure}[t]
\begin{center}
\includegraphics[width=4.5in]{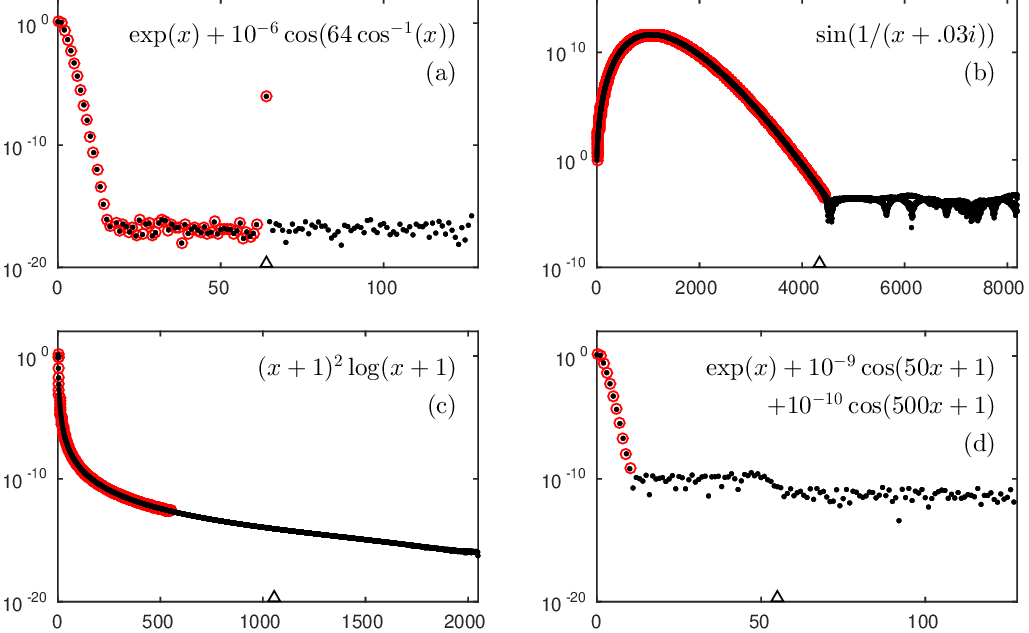}
\caption{\label{fig8}\small(a) {\tt sampleTest} catches the $T_{64}$ component.
(b) Some series grow exponentially before decaying.
(c) When convergence is slow, Chebfun balances cost and
benefit and chops before reaching the target tolerance.
(d) Here ${\tt plateauPoint} = 56$ but ${\tt cutoff} = 11$.}
\end{center}
\end{figure}

Example \ref{fig7}a repeats $\exp(x)$, but with a little noise added
in the form of a term $10^{-14} \cos(999 x)$, raising the
level of the plateau.  Here the flexibility
of the constructor leads to a chop one term earlier than before
(algorithmic aim~6).
Examples \ref{fig7}b--\ref{fig7}d explore the engineering choices involved in
the design of the plateau detector.  In Example \ref{fig7}b, a short
plateau at level $10^{-10}$ is not long enough to count, and
the constructor insists on going further.  The longer
quite flat plateau at the same level in Example \ref{fig7}c, on the other hand, is
accepted for chopping, a consequence of
algorithmic aim~9.  Example \ref{fig7}d shows that raising
the magnitude by one further factor of $10$ is enough for Chebfun
once again not to accept it as a plateau (algorithmic aim~5).

Example \ref{fig8}a illustrates how {\tt sampleTest} can catch a component
that otherwise might have been missed because of aliasing.
Example \ref{fig8}b shows that the scale of a function may be determined by
Chebyshev coefficients far from the beginning of the series.  This 
example also illustrates Chebfun's treatment of complex functions.
In Example \ref{fig8}c the series converges very slowly, and Step~3 of
{\tt standardChop} leads to a reduction of {\tt cutoff}
well below {\tt plateauPoint}, the rationale being that the last
500 coefficients have not brought enough improvement to be
worth the cost.  (To force higher accuracy,
one could specify the length of the 
chebfun explicitly.)  Example \ref{fig8}d also has
${\tt cutoff} \ll {\tt plateauPoint}$, but here it is
a matter of avoiding a chopping point that would have been 
``obviously wrong'' (algorithmic aim 11).

\begin{figure}[t]
\begin{center}
\includegraphics[width=4.5in]{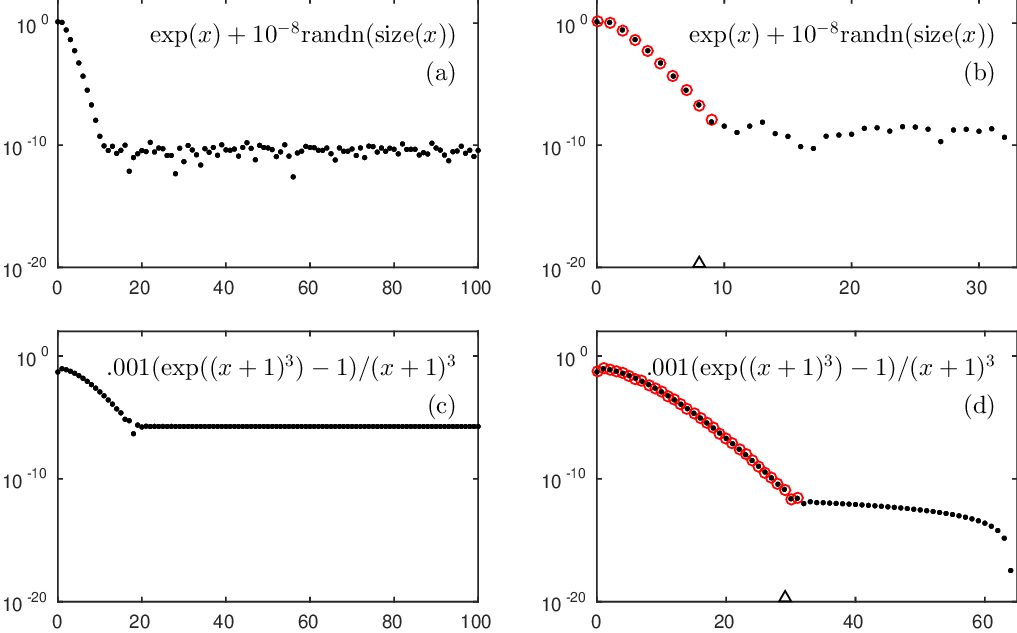}
\caption{\label{fig9}\small(a) Explicit noise prevents convergence.
(b) Convergence achieved with ${\tt eps} =10^{-8}$.
(c) Noise from cancellation error prevents convergence.
(d) Convergence achieved with ${\tt eps} =10^{-10}$.}
\end{center}
\end{figure}

The examples of Figure~\ref{fig9} explore the treatment of noise
(algorithmic aim~8).
Example \ref{fig9}a involves a function with too much noise;
{\tt standardChop} is never happy on any grid, and Chebfun issues a
warning message after attempting
the default maximum grid size of $2^{16}+1$.  Example \ref{fig9}b
shows the different result achieved if {\tt eps} is
increased to $10^{-8}$.  Examples \ref{fig9}c and \ref{fig9}d are
a similar pair, except now, the noise is introduced not explicitly
by {\tt randn} but implicitly via cancellation error
in the computation of a function
with a removable singularity at $x=-1$.

\section{Standard, strict, and loose constructors}
The name {\tt standardChop} comes from an idea of long standing
in the Chebfun project, that a user may
wish to apply different chopping strategies in different contexts.
Sometimes one might want to insist on finding coefficients
beneath a specified tolerance: a {\em ``strict''} chopping rule.
Sometimes one might wish to seek a plateau of noise or
rounding errors without any a priori notion whatever about how
high this might be: a {\em ``loose''} chopping rule.  Our presumption
has been that most of the time, one would wish to operate in
an in-between, {\em ``standard''} mode, where the algorithm
has an a priori notion of {\tt eps} but applies it
with some flexibility.  In {\tt standardChop,} as we have
described, the series ultimately accepted will normally stop
somewhere between ${\tt eps}^{2/3}$ and ${\tt eps}$.

Before the release of Version 5.3, Chebfun offered
strict and standard (=``classic'') construction modes,
with the latter as the default.  The strict mode got almost
no use, however, and a loose mode was not implemented.
The new {\tt standardChop} has now replaced ``classic'' as the default.

\section{Computation, construction, and simplification}
Up to now, we have written as if all Chebfun does is construct
polynomial approximations from strings or anonymous functions.
But if this were the case, there would be no need for it!---the
job could be done better with symbolic computing.
The point of a system like this is not just construction of
functions from scratch, but numerical computation with 
functions, which necessitates ongoing use of the construction
process at each step.

Most of the time, Chebfun follows this prescription exactly.
If a chebfun is operated on by an operation like {\tt exp(f)}, or if
two chebfuns are combined by an operation like {\tt f./g}, the
usual process consists of using the existing representations
of {\tt f} and {\tt g} to compute sample values on grids of
size $17, 33, \dots,$ and apply {\tt standardChop} as usual to
construct the new result.
This description applies both for classic nonperiodic chebfuns,
represented by Chebyshev series, and the more recently introduced
periodic chebfuns, as discussed in the next section.

Some mathematical operations, on the other hand,
can (or must) be carried out without
going through this construction process.  
The simplest example is the unary minus operation: to construct {\tt -f}
from {\tt f}, we negate the coefficients rather than sampling {\tt f}
anew.  Besides taking a little time,
resampling would sometimes change the representation
slightly since it would
introduce new rounding errors.  For example,
the chebfun {\tt f = chebfun(@(x) sin(1./(x+.03i))} in Figure \ref{fig8}b
has length 4441, and so does {\tt -f}, but if we construct the negative
via {\tt chebfun(@(x) -f(x))} the rounding error
plateau is higher and the length reduces to $4266$.\footnote{What if
we keep reconstructing over and over again?  The lengths come out
as $4262, 4262, 4261, 4261, 4261, 4261,\dots.$  Thus chebfun construction
is approximately, but not exactly, a projection.}

Some operations are not carried out by calling the constructor,
but lead to Chebyshev series whose tails may be below the noise level.
One instance of this is {\tt cumsum}, which computes the integral of a chebfun
by manipulating its Chebyshev coefficients.  The amplitudes of
the higher coefficients will be reduced in the process and can
be trimmed without loss of overall accuracy.  This is intuitively natural
since integration is a smoothing operation and thus can be expected
to result in shorter series.  Chebfun has a
code {\tt simplify} that effects this trimming by calling
{\tt standardChop} after first prolonging the length of the Chebyshev series 
by about $25\%$ to ensure it is long enough for {\tt standardChop}
to work properly.
For example,
{\tt f = chebfun(@(x) log(1.1-x))} has length 75, so mathematically, one might
expect {\tt cumsum(f)} to have length 76.  Since
{\tt cumsum} calls {\tt simplify}, however, the length is
actually~70.  

Multiplication of two chebfuns, the operation {\tt f.*g}, is also
currently carried out in Chebfun by manipulating coefficients, then
calling {\tt simplify}.
Addition and subtraction are currently carried out
by manipulating coefficients without calling {\tt simplify}, but this
is likely to change.

Another example of a Chebfun algorithm that bypasses the
construction process is {\tt conv,} which computes the convolution 
of two functions.  This code uses a fast algorithm introduced
by Hale and Townsend that converts Chebyshev series to Legendre
series for the convolution, then converts back again and calls
{\tt simplify}~\cite{haletown}.
Construction via pointwise sampling would be much more expensive.

\section{Variant construction processes}
We have presented a general construction process whereby Chebfun
samples a function on successively finer grids, computes Chebyshev
coefficients, and calls {\tt standardChop} to decide if convergence
has been achieved.  We now describe various parts of Chebfun
in which this algorithm is applied in specialized ways.
Before the introduction of Version 5.3, these operations were
carried out by different alogorithms, and their replacement by
{\tt standardChop} has been a significant simplification
(algorithmic aim~4 of Section~4).

\subsection{Chebfuns with several pieces}
A chebfun may consist of a concatenation of several
smooth pieces, called {\em funs,} each with its own
Chebyshev series representation.  Breakpoints between
funs are introduced by operations known to break
smoothness at certain points, such as {\tt abs(f)} at points
where {\tt f} passes through zero, and they can also be
determined by a fast edge-detection algorithm~\cite{ppt}.

The complication that arises in construction of chebfuns
with several pieces is that
a chebfun can normally be evaluated to a certain accuracy relative
to its global scale over the whole interval of definition---not
the local scale of an individual fun.  To achieve
the necessary effect, {\tt standardChop} is called with
its input parameter {\tt tol} set not to the Chebfun tolerance
{\tt eps} as usual, but to {\tt eps*vscaleGlobal/vscaleLocal},
where {\tt vscaleGlobal} and {\tt vscaleLocal} are estimates
of the scales of the global chebfun and the local fun, respectively.

\subsection{Periodic functions and Fourier series}
Beginning with Version 5.0 in 2014, Chebfun has had a capability
of working with periodic functions represented by Fourier series
in addition to the usual nonperiodic functions represented by
Chebyshev series~\cite{trigfun}. This adds efficiency and accuracy for periodic
problems and is particularly attractive in eliminating
discontinuities of periodic functions or their derivatives
at boundary points.  A~periodic
chebfun is informally called a ``trigfun''.

\begin{figure}
\begin{center}
\includegraphics[scale=.62]{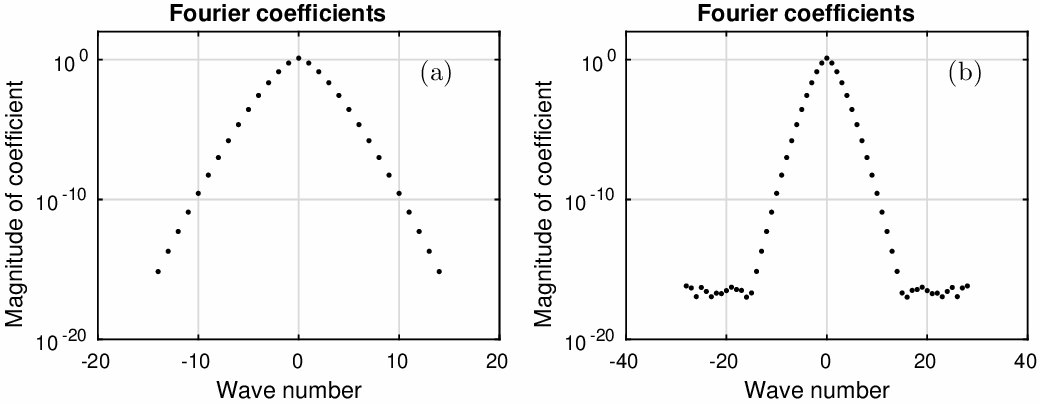}
\end{center}
\vskip -1em
\caption{\label{doublelengthtrig}\small Repetition of
Figure~\ref{doublelength} for the function
$f(x) = \exp(\sin(\pi x))$, but now with Chebfun called
with the {\tt 'trig'} flag, producing a Fourier rather
than Chebyshev representation.}
\end{figure}

To construct a trigfun, which occurs when the user specifies
the flag \verb|'trig'|, Chebfun 
samples the function on periodic equispaced
grids of $16, 32, \dots$ points and constructs trigonometric series for
the corresponding interpolating trigonometric polynomials.  These
series are of degrees $8, 16, \dots$ in the usual terminology, where
a series of degree 8 for $t\in [-\pi,\pi]$,
for example, consists of coefficients 
$a_{-8},\dots , a_8$ multiplying exponentials $e^{-8it},\dots, 
e^{8it}$.  The decision of when a series is ``happy''
and where to chop it is made by {\tt standardChop}
applied to the sequence
$|c_0|$, $(|c_1|+|c_{-1}|)/2$ (repeated twice),
$(|c_2|+|c_{-2}|)/2$ (repeated twice), and so on.
The reason for this duplication of
each coefficient is so that Fourier series will
be treated by essentially the same parameters as Chebyshev series,
with 17 values being the minimal number for happiness.
An example of Fourier coefficients of a trigfun is shown in
Figure~\ref{doublelengthtrig}, which revisits the function
$f(x) = \exp(\sin(\pi x))$ of Figure~\ref{doublelength} now in
Fourier mode.

\subsection{ODE boundary-value problems}
One of the most important capabilities of Chebfun for users
is the solution of ordinary differential equations (ODEs), 
both boundary-value problems (BVPs) and initial-value 
problems (IVPs), which can be linear or nonlinear.  

To solve a linear BVP, Chebfun discretizes the problem on
Chebyshev grids of sizes approximately
$33, 65, 129, 257, 513, 725, 1025, 1449,\dots$ and checks for happiness
on each grid~\cite{dbt,dh}.
Apart from the modified grid sequence, which is based on half-integer
as well as integer powers of~2, this
process differs from standard construction
in two ways.  One is that this is not just a matter of ``sampling''
a fixed function, since the value at a grid point such as $0$, for example,
will change as the gridding is refined.  In Chebfun terminology,
this means BVPs are constructed in a ``resampling'' mode, with
function values already obtained
on a coarse grid recomputed on each finer grid.
The other difference is that for solving BVPs, {\tt standardChop}
is called with ${\tt tol} = {\tt bvpTol}$, where the parameter
{\tt bvpTol} is by default set to $5\times 10^{-13}$
rather than machine epsilon. One
reason for this is that solution of BVPs on fine grids is expensive,
with $O(n^3)$ complexity on a grid of $n$ points, so pushing
to full machine precision may be slow.  In addition, 
the matrices involved in the solution process are often ill-conditioned,
so setting ${\tt tol} = {\tt eps}$ would sometimes be problematic.

\begin{figure}
\begin{center}
\includegraphics[width=3.2in]{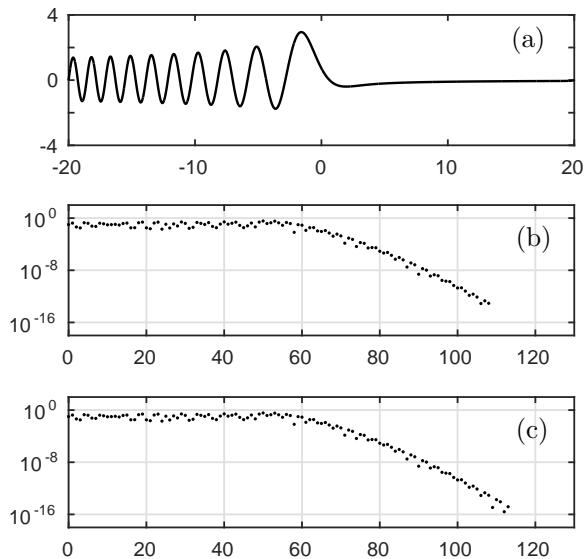}
\end{center}
\vskip -1.5em
\caption{\small\label{figode}(a) Chebfun solution to
$u'' - xu = 1$ for $x\in [-20,20]$.
(b) Chebyshev coefficients $|a_k|$ for standard construction with
the default tolerance 
${\tt bvpTol} = 5\times 10^{-13}$. 
(c) More accuracy achieved by
tightening the tolerance to 
${\tt bvpTol} = 10^{-16}$.}
\end{figure}

Figure~\ref{figode} illustrates the Chebfun solution of
$u'' - xu = 1$ for $x\in [-20,20]$ with $u(-20)=u(20)=0$ and
the corresponding amplitudes of Chebyshev coefficients.  In this
unproblematic
case the default tolerance ${\tt tol} = 5\times 10^{-13}$ does
not achieve the maximum accuracy possible, and a second
plot of coefficients is shown for a
solution with ${\tt bvpTol} = 10^{-16}$.

If a BVP is nonlinear, Chebfun embeds the construction process
just described in a Newton or damped-Newton iteration, with
the necessary derivatives formulated in a continuous setting
as Fr\'echet derivative operators constructed by automatic
differentiation~\cite{bd}.  The correction functions
will eventually become very small,
and it is necessary for the Chebyshev series involved in
their construction to be judged
relative to the scale of the overall function, not that of the
correction.  Accordingly, {\tt standardChop} is called with
its value of {\tt tol} increased by a factor analogous to
{\tt vscaleGlobal/vscaleLocal}, as described in \S 8.1.
The Newton iteration is exited when its estimated error falls
below $200\times {\tt bvpTol}$.

\subsection{ODE initial-value problems}
IVPs in Chebfun are solved differently from BVPs, by
marching rather than a global spectral
discretization~\cite{birk}.  This
solution process has nothing a priori to do with Chebyshev grids,
and it is carried out by standard Matlab ODE codes: {\tt ode113} by
default, which can be changed e.g.\ to {\tt ode15s} for
a stiff problem.  As with BVPs, Chebfun aims by default for about
12 digits rather than 16.  To be precise, {\tt ode113} by
default is called with
{\tt absTol = 1e5*macheps} and {\tt relTol = 1e2*macheps} (see
footnote~2).  The computed
function that results is then converted to a chebfun by a call
to {\tt standardChop} with {\tt tol} set to the maximum of
{\tt relTol} and {\tt absTol/vscale}.

\subsection{Quasimatrices}
A quasimatrix is a chebfun with more than one column, that is,
a collection of several functions defined on the same interval
$[a,b]$.  By default, Chebfun constructs each column of 
a quasimatrix independently, calling {\tt standardChop}
with ${\tt tol}={\tt eps}$ as usual.  For some applications, however,
it is appropriate for the columns to be constructed relative
to a single global scale, for example if they correspond to pieces of
a decomposition of a single function.
For such applications a user can specify the
flag \verb|'globaltol'|, and then {\tt standardChop} will be
called with {\tt tol} appropriately adjusted for the various
columns as described in \S 8.1 and 8.3.  The next subsection
gives an example.

\begin{figure}
\begin{center}
\includegraphics[width=3.1in]{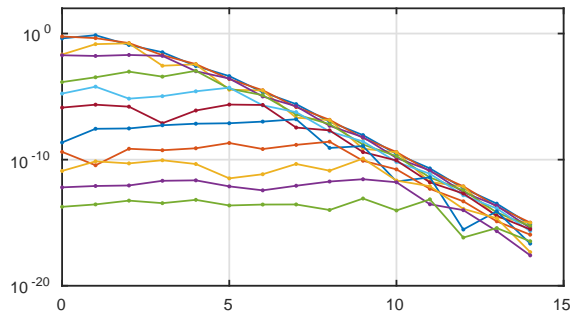}
\end{center}
\vskip -1.6em
\caption{\small\label{figquasi}Chebyshev coefficients of the
rows of the chebfun2 of $f(x,y) = \cos(xy+1)$, corresponding to the
$x$ part of the bivariate representation.
The curves for the 12 rows begin at different heights but end at approximately
the same height, reflecting calls to 
{\tt standardChop} with the {\tt 'globaltol'} flag set
so that tolerances are adjusted to a global scale.}
\end{figure}

\subsection{Multiple dimensions: Chebfun2 and Chebfun3}
Finally we mention that Chebfun can also compute with smooth functions 
in two dimensions, since the release of Chebfun2 in
2012~\cite{chebfun2}, and soon in 3D with the upcoming
release of Chebfun3~\cite{chebfun3}.  In both cases, functions are represented
by low-rank approximations constructed from outer products
of 1D functions, which in turn are represented by the usual Chebfun
Chebyshev series, or Fourier series in periodic directions.
These series are constructed with {\tt standardChop} in the
\verb|'globaltol'| mode described above.  To illustrate,
Figure~\ref{figquasi}
plots the Chebyshev coefficients of the 12 rows representing
the $x$-dependence of the Chebfun2
representation of $f(x,y) = \cos(xy+1)$.

\section{Conclusions}
Although Chebfun has been chopping Chebyshev series for
more than a decade, and more recently also Fourier series,
the details have been inconsistent and ad hoc until lately.
This paper has described the new algorithm 
{\tt standardChop}, which unifies these processes
with a clear structure.
Other related projects for numerical computation with functions
face the same challenge of ``rounding'' of functions.
A list of Chebfun-related projects can be found under
the ``About'' tab at {\tt www.chebfun.org,} currently including
ApproxFun.jl~\cite{approxfun}, pychebfun,
Fourfun~\cite{fourfun}, CHEBINT~\cite{chebint}, PaCAL~\cite{pacal},
sincfun~\cite{sincfun},
a collection of LISP codes by Fateman~\cite{fateman},
libchebfun, and rktoolbox~\cite{rktoolbox}.

The introduction of {\tt standardChop} has
made the foundations of the Chebfun project more secure.  At the
same time, we reiterate that
we make no claim that this algorithm is optimal, or even
represents the only reasonable approach to this problem. Just as
it took decades for floating-point arithmetic of real numbers to reach
a reasonably settled state with the introduction of the IEEE
Floating Point Standard in the mid-1980s~\cite{overton},
perhaps it will take a long time
for consensus to emerge as to the best ways to realize the
analogue of floating-point arithmetic for numerical computation
with functions.

\section*{Acknowledgments}
We have benefitted from extensive discussions with 
others in the Chebfun team including Anthony Austin,
\'Asgeir Birkisson,
Toby Driscoll, Nick Hale, Behnam Hashemi, Mohsin
Javed, Hadrien Montanelli, Mikael Slevinsky,
Alex Townsend, Grady Wright, and Kuan Xu.
We also acknowledge insightful suggestions from Folkmar
Bornemann of TU Munich.
This work was supported by the European Research Council
under the European Union's Seventh Framework Programme
(FP7/2007--2013)/ERC grant agreement no.\ 291068.\ \ The
views expressed in this article are not those of the ERC
or the European Commission, and the European Union
is not liable for any use that may be made of the information
contained here.

\section*{Appendix: listing of {\tt standardChop.m}}

{\footnotesize
\begin{verbatim}
function cutoff = standardChop(coeffs, tol)
%STANDARDCHOP  A sequence chopping rule of "standard" (as opposed to "loose" or
% "strict") type, that is, with an input tolerance TOL that is applied with some
% flexibility.  This code is used in all parts of Chebfun that make chopping
% decisions, including chebfun construction (CHEBTECH, TRIGTECH), solution of
% ODE BVPs (SOLVEBVP), solution of ODE IVPs (ODESOL), simplification of chebfuns
% (SIMPLIFY), and Chebfun2.  See J. L. Aurentz and L. N. Trefethen, "Chopping a
% Chebyshev series", arXiv, December 2015.
%
% Input:
%
% COEFFS  A nonempty row or column vector of real or complex numbers
%         which typically will be Chebyshev or Fourier coefficients.
%
% TOL     A number in (0,1) representing a target relative accuracy.
%         TOL will typically will be set to the Chebfun EPS parameter,
%         sometimes multiplied by a factor such as vglobal/vlocal in
%         construction of local pieces of global chebfuns.
%         Default value: machine epsilon (MATLAB EPS).
%
% Output:
%
% CUTOFF  A positive integer.
%         If CUTOFF == length(COEFFS), then we are "not happy":
%         a satisfactory chopping point has not been found.
%         If CUTOFF < length(COEFFS), we are "happy" and CUTOFF
%         represents the last index of COEFFS that should be retained.
%
% Examples:
%
% coeffs = 10.^-(1:50); random = cos((1:50).^2);
% standardChop(coeffs) % = 18
% standardChop(coeffs + 1e-16*random) % = 15
% standardChop(coeffs + 1e-13*random) % = 13
% standardChop(coeffs + 1e-10*random) % = 50
% standardChop(coeffs + 1e-10*random, 1e-10) % = 10
 
% Jared Aurentz and Nick Trefethen, July 2015.
%
% Copyright 2015 by The University of Oxford and The Chebfun Developers. 
% See http://www.chebfun.org/ for Chebfun information.

% STANDARDCHOP normally chops COEFFS at a point beyond which it is smaller than
% TOL^(2/3).  COEFFS will never be chopped unless it is of length at least 17 and
% falls at least below TOL^(1/3).  It will always be chopped if it has a long
% enough final segment below TOL, and the final entry COEFFS(CUTOFF) will never
% be smaller than TOL^(7/6).  All these statements are relative to
% MAX(ABS(COEFFS)) and assume CUTOFF > 1.  These parameters result from
% extensive experimentation involving functions such as those presented in
% the paper cited above.  They are not derived from first principles and
% there is no claim that they are optimal.

% Set default if fewer than 2 inputs are supplied: 
if ( nargin < 2 )
    p = chebfunpref;
    tol = p.eps;
end

% Check magnitude of TOL:
if ( tol >= 1 ) 
    cutoff = 1;
    return
end

% Make sure COEFFS has length at least 17:
n = length(coeffs);
cutoff = n;
if ( n < 17 )
    return
end
  
% Step 1: Convert COEFFS to a new monotonically nonincreasing
%         vector ENVELOPE normalized to begin with the value 1.

b = abs(coeffs);
m = b(end)*ones(n, 1);
for j = n-1:-1:1
    m(j) = max(b(j), m(j+1));
end   
if ( m(1) == 0 )
    cutoff = 1;
    return
end
envelope = m/m(1);

% Step 2: Scan ENVELOPE for a value PLATEAUPOINT, the first point J-1, if any,
% that is followed by a plateau.  A plateau is a stretch of coefficients
% ENVELOPE(J),...,ENVELOPE(J2), J2 = round(1.25*J+5) <= N, with the property
% that ENVELOPE(J2)/ENVELOPE(J) > R.  The number R ranges from R = 0 if
% ENVELOPE(J) = TOL up to R = 1 if ENVELOPE(J) = TOL^(2/3).  Thus a potential
% plateau whose starting value is ENVELOPE(J) ~ TOL^(2/3) has to be perfectly
% flat to count, whereas with ENVELOPE(J) ~ TOL it doesn't have to be flat at
% all.  If a plateau point is found, then we know we are going to chop the
% vector, but the precise chopping point CUTOFF still remains to be determined
% in Step 3.

for j = 2:n
    j2 = round(1.25*j + 5); 
    if ( j2 > n )
        % there is no plateau: exit
        return
    end      
    e1 = envelope(j);
    e2 = envelope(j2);
    r = 3*(1 - log(e1)/log(tol));
    plateau = (e1 == 0) | (e2/e1 > r);
    if ( plateau )
        % a plateau has been found: go to Step 3
        plateauPoint = j - 1;
        break
    end
end

% Step 3: fix CUTOFF at a point where ENVELOPE, plus a linear function
% included to bias the result towards the left end, is minimal.
%
% Some explanation is needed here.  One might imagine that if a plateau is
% found, then one should simply set CUTOFF = PLATEAUPOINT and be done, without
% the need for a Step 3. However, sometimes CUTOFF should be smaller or larger
% than PLATEAUPOINT, and that is what Step 3 achieves.
%
% CUTOFF should be smaller than PLATEAUPOINT if the last few coefficients made
% negligible improvement but just managed to bring the vector ENVELOPE below the
% level TOL^(2/3), above which no plateau will ever be detected.  This part of
% the code is important for avoiding situations where a coefficient vector is
% chopped at a point that looks "obviously wrong" with PLOTCOEFFS.
%
% CUTOFF should be larger than PLATEAUPOINT if, although a plateau has been
% found, one can nevertheless reduce the amplitude of the coefficients a good
% deal further by taking more of them.  This will happen most often when a
% plateau is detected at an amplitude close to TOL, because in this case, the
% "plateau" need not be very flat.  This part of the code is important to
% getting an extra digit or two beyond the minimal prescribed accuracy when it
% is easy to do so.

if ( envelope(plateauPoint) == 0 )
    cutoff = plateauPoint;
else
    j3 = sum(envelope >= tol^(7/6));
    if ( j3 < j2 )
        j2 = j3 + 1;
        envelope(j2) = tol^(7/6);
    end
    cc = log10(envelope(1:j2));
    cc = cc(:);
    cc = cc + linspace(0, (-1/3)*log10(tol), j2)';
    [~, d] = min(cc);
    cutoff = max(d - 1, 1);
end

end
\end{verbatim}
\par}

\end{document}